\documentclass[fleqn]{mat01}
\usepackage{times,mathtimy,amssymb}
\begin{document}

\setcounter{page}{181}
\firstpage{181}

\newtheorem{theore}{Theorem}
\renewcommand\thetheore{\arabic{section}.\arabic{theore}}
\newtheorem{theor}[theore]{\bf Theorem}
\newtheorem{lem}[theore]{Lemma}
\newtheorem{exampl}[theore]{Example}

\newtheorem{theorr}{Theorem}
\renewcommand\thetheorr{\arabic{theorr}}
\newtheorem{therr}[theorr]{\bf Theorem}
\newtheorem{coro}[theorr]{\rm COROLLARY}
\newtheorem{propo}[theorr]{\rm PROPOSITION}
\newtheorem{rem}[theorr]{Remark}
\newtheorem{case}[theorr]{\it Case}
\def\lemm{\trivlist\item[\hskip\labelsep{{\it Lemma.}}]}
\def\remm{\trivlist\item[\hskip\labelsep{{\it Remark.}}]}
\def\definit{\trivlist\item[\hskip\labelsep{{\rm DEFINITION}}]}

\font\wwww=tibi at 13.5pt
\def\n{\mbox{\wwww{n}}}

\title{Geometry of good sets in $\n$-fold Cartesian product}

\markboth{A~K{\l}opotowski, M~G~Nadkarni and K~P~S~Bhaskara Rao}{Geometry of good sets in $n$-fold Cartesian product}

\author{A~K{\L}OPOTOWSKI$^{1}$, M~G~NADKARNI$^{2,3}$ and K~P~S~BHASKARA~RAO$^{4}$}

\address{$^{1}$Universit{\'e} Paris XIII, Institut Galil{\'e}e, 93430
Villetaneuse Cedex, France\\
\noindent $^{2}$Institute of Mathematical Sciences, C.I.T. Campus,
Chennai 600 113, India\\
\noindent $^{3}$Chennai Mathematical Institute, Chennai 600 117, India\\
\noindent $^{4}$Department of Mathematics, Southwestern College,
Winfield, KS 76156, USA\\
\noindent E-mail: klopot@math.univ-paris13.fr;
nadkarni@math.mu.ac.in; kpsbrao@hotmail.com}

\volume{114}

\mon{May}

\parts{2}

\Date{MS received 16 November 2003; revised 3 February 2004}

\begin{abstract}
We propose here a multidimensional generalisation of the notion of link
introduced in our previous papers and we discuss some consequences for
simplicial measures and sums of function algebras.
\end{abstract}


\maketitle

\setcounter{section}{-1}
\section{Introduction}

Let $X_1, X_2, \ldots, X_n$ be non-empty sets and let $\Omega =
X_1\times X_2 \times \cdots \times X_n$ be their Cartesian product. For
each $i, 1\leq i \leq n$, $\Pi_i$ will denote the canonical projection
of $\Omega$ onto $X_i$. A subset $S \subset \Omega$ is said to be {\it
good} if every complex valued function $f$ on $S$ is of the\break form:
\begin{align*}
f(x_1, x_2, \ldots, x_n) &= u_1(x_1) + u_2(x_2) + \cdots + u_n(x_n), (x_1, x_2, \ldots, x_n) \in S,
\end{align*}
for suitable functions $u_1, u_2, \ldots, u_n$ on $X_1, X_2, \ldots,
X_n$ respectively.

A necessary and sufficient condition for a subset $S$ of $X_1\times
X_2\times\cdots\times X_n$ to be good was derived in our paper \cite{7}
and some consequences for simplicial measures and sums of algebras were
discussed. For $n=2$ these questions are well-discussed in
\cite{1,2,3,5,6,7,10,11,12,13,14,17}. The notion of a link or path
between two points plays a crucial role in all these papers. For $n >2$
a natural notion of link between two points of $S$ was so far not
available, a difficulty mentioned on p.~82 and 84 of \cite{7}. So
natural analogues of results for $n=2$ were not available for the case
$n >2$. This paper attempts to remove this difficulty. Here we define,
for $n \geq 2$, what we call full sets in terms of which a notion of
geodesic between two points of a good set is formulated. This allows us
to prove some results on simplicial measure and sums of algebras in
terms of geodesics in analogy with the case $n = 2$. For $n = 2$ a
geodesic between two points is a link as defined in \cite{3}, and for $n
> 2$ a geodesic has nearly all the properties of this object. For
question concerning sums of algebras for $n >2$ we refer to the papers
\cite{18,19} where the notions of uniformly separating families and
uniformly measure separating families are introduced and applied both
for questions of sums of algebras and in dimension theory, and to paper
\cite{16}.

\pagebreak

\section{Examples}

\begin{enumerate}
\leftskip .5pc
\item A singleton subset of $\Omega$ is always a good set. Also any
subset of $\Omega$ no two points of which have a coordinate in common
is a good set.

\item The subset $S = \{(0,0), (1,0),(0,1)\}$ of $\{0,1\}\times \{0,1\}$
is a good set. For let $f$ be any function on $S$ and let $u_1(0)$ be
given an arbitrary value, say $c$, and define $u_2(0) = f(0,0) - c$.
With $u_2(0)$ thus defined, we write $u_1(1) = f(1,0) - u_2(0)$. Finally
we get $u_2(1) = f(1,1) -u_1(1)$. Clearly $u_1 + u_2 = f$ on $S$. Note
that once $u(0)$ is fixed, the solution is unique.

\item Let $S\subset X_1\times X_2$. Say that two points $(x,y), (z,w)$
in $S$ are linked if there is a finite sequence $(x_1,y_1),(x_2,y_2),
\ldots, (x_n,y_n)$ in $S$ such that (i) $(x_1,y_1)=(x,y), (x_n,y_n) =
(z,w)$, (ii) for each $i, 1\leq i \leq n - 1$, exactly one of the two
inequalities holds $x_i \neq x_{i+1}, y_i \neq y_{i+1}$, (iii) if for
any $i$, $x_i \neq x_{i+1}$ then $x_{i+1} = x_{i+2}$ and if $y_i \neq
y_{i+1}$ then $y_{i+1} = y_{i+2}$, $1 \leq i \leq n - 2$. If $(x,y)$ and
$(z,w)$ are linked we write $(x,y)L(z,w)$ and observe that $L$ is an
equivalence relation. If there is only one link between two points
$(x,y)$ and $(z,w) \in S$, then we say that $(x,y)$ and $(z,w)$ are
uniquely linked. We note that $S$ is good if and only if any two linked
points in $S$ are uniquely linked. If $S$ is good and $C$ is a set which
meets each equivalence class of $L$ in exactly one point, then the
solution of $u_1(x_1) + u_2(x_2) = f(x_1,x_2)$ is unique once we
prescribe the values of $u_1$ on $\Pi_1C$ (see \cite{3}).

\item The set $\{(0,0,0), (1,0,0), (1,1,0), (1,1,1), (2,1,1), (2,2,1),
\ldots\}$ where starting at $(0,0,0)$ one moves one unit at a time,
first along the $x$-axis, then along the $y$-axis and then along the
$z$-axis and continuing similarly with the next movement along the
$x$-axis, is a good set. For any $f$ on this set, the solution of
$u_1(x_1) + u_2(x_2) + u_3(x_3) = f(x_1,x_2, x_3)$ is unique once we
prescribe the values of $u_1(0)$ and\break $u_2(0)$.

\item $S = \{(0,0,0), (1, 0, 0), (0,1,0), (0,0,1)\}$ is a good set in
$\{0,1\}^3$ while the set $S\cup \{(1,1,1)\}$ is not a good set.

\item $S = \{(1,1,0), (1,0,1), (0,1,1), (0,0,0)\}$ is a good set in $\{
0,1\}^3$. This example is different from example 4 in that no two
elements of $S$ differ from each other in only one coordinate, yet for
any $f$, the solution of $u_1(x_1) + u_2(x_2) + u_3(x_3) = f(x_1,x_2,
x_3)$ is unique once we prescribe the values of $u_1(0)$ and $u_2(0)$.

\item $\{(1,2,3),(4,5,6), (7,8,9), (1,5,9)\}$ is a good set. For a given
$f$ on $S$, the equation $u_1(x_1) + u_2(x_2) + u_3(x_3) = f(x_1,x_2,
x_3), (x_1,x_2,x_3)\in S$ gives four linear equations in nine variables.
If we fix the values of some suitable five variables, then the solution
is unique, but not any choice of five variables would do.

\item Let $a_i \in X_i$, $i = 1,2,3$. Then 
\begin{equation*}
\hskip -.5cm S = X_1\times \{a_2\}\times \{a_3\}\cup \{a_1\}\times X_2
\times \{a_3\}\cup \{a_1\}\times \{a_2\}\times X_3
\end{equation*}
is a good set in $X_1\times X_2 \times X_3$.

\item The embedding of the $n$-dimensional unit cube $E^n$ into
${\mathbb{R}}^{2n+1}$ obtained in Kolmogorov's solution of Hilbert's
thirteenth problem \cite{8} is a good set.

\item If $S$ is a good set in $X_1\times X_2$ and $(x_0,y_0) \in S$ then
$U, V$ which satisfy $u(x) + v(y) = 1_{\{(x_0,y_0)\}}(x,y), ~(x,y)\in
S$, $u(x_0) = 0$ are necessarily bounded in absolute value by 1.
However, this can fail if $n >2$ as the following example, obtained
jointly with Gowri Navada, shows: Consider the set $\{(x_0,y_0,z_0)$,
$(x_1, y_0,z_0)$, $(x_0, y_1 ,z_0)$, $ (x_1, y_1,z_1)$, $ (x_2, y_0,z_1)
$, $(x_0, y_2,z_1)$, $(x_2, y_2,z_2)$,$\ldots$, $(x_n,y_n,z_n)$,
$(x_{n+1},y_{0}, z_n)$, $(x_{0},y_{n+1},z_n)$, $(x_{n+1},y_{n+1},
z_{n+1})$, $\ldots\}$ in $X\times Y \times Z$, where $X$, $Y$, $Z$ are
infinite sets. This is a good set since each point admits a coordinate
which does not appear as a coordinate of any of the points preceding
it. Further it is easily seen that the solution $U, V, W$ of 
\begin{equation*}
\hskip -.5cm u(x) + v(y) + w(z) = 1_{\{(x_0,y_0,z_0)\}}(x,y,z),\quad
(x,y,z) \in S,
\end{equation*}
satisfying $u(x_0) = 0, v(y_0) = 0$, is given by, $W(z_0) = 1$ and for
$n > 0$, $U(x_n) = V(y_n) = -2^{n-1}, W(z_{n}) = 2^n$.
\end{enumerate}

\section{Characterisation of good sets; consequences}

Given any finitely many symbols $t_1, t_2, \ldots, t_k$ with {\it
repetitions} allowed and given any finitely many integers $n_1, n_2,
\ldots, n_k$, we say that the formal sum $n_1t_1 + n_2t_2 + \cdots +
n_kt_k$ vanishes if for every $t_j$ the sum of the coefficients of $t_j$
vanishes.

\begin{definit}$\left.\right.$\vspace{.5pc}

\noindent {\rm An element $(x_1, x_2, \ldots, x_n)$ of $\Omega$ will be
denoted by $\vec x$. A non-empty finite subset $ L = \{{\vec {x}}_1,
{\vec{x}}_2, \ldots, {\vec{x_k}}\}$ of $\Omega$ is called a {\it loop}
if there exist non-zero integers $n_1, n_2, \ldots, n_k$ such that the
sum $\sum_{i=1}^kn_i{\vec{x}}_i$ vanishes in the sense that the formal
sum vanishes coordinatewise, and no strictly smaller non-empty subset
of $L$ has this property.}
\end{definit}\vspace{.5pc}

We have {\it $S \subset \Omega$ is good if and only if there are no
loops in $S$.} This characterisation of a good set, proved in \cite{7},
implies:
\begin{enumerate}
\renewcommand\labelenumi{(\arabic{enumi})}
\leftskip .1pc
\item $S$ is good if and only if every finite subset of $S$ is good,

\item union of any directed family of good sets is a good set, where a
family of sets is said to be directed if given any two sets in the
family there is a third set in the family which includes both. In
particular, any union of a linearly ordered (under inclusion) system of
good sets is a good set,

\item in view of (2), by Zorn's lemma, we conclude that every good set
is contained in a maximal good set, where a good subset in $\Omega$ is
said to be maximal if it is not contained in a strictly larger good
subset of $\Omega$.\vspace{-1pc}
\end{enumerate}
Note that if $S\subset \Omega$ is maximal then, for each $i$, $\Pi_i S =
X_i$, for if $X_i - \Pi_i S$ is non-empty for some $i$, and if ${\vec
{x}}\in \Omega$ has $i$th coordinate not in $\Pi_i S$, then $S
\cup\{\vec{x}\}$ is a good set bigger than $S$.

\section{Full sets}

The following refined notion of maximal set, called full set, will be
crucial for our discussion.

\begin{definit}$\left.\right.$\vspace{.5pc}

\noindent {\rm A subset $S$ of $\Omega$ is said to be full if $S$ is a
maximal good set in $\Pi_1 S\times \Pi_2 S\times\cdots\times \Pi_nS$.}
\end{definit}\vspace{.5pc}

Clearly every good set $S$ is contained in a full good set $S'$ such
that the canonical projections of $S$ and $S'$ on the coordinate spaces
coincide.

\begin{therr}[\!] Let $S\subset \Omega$ be given. Assume that there
exist $x_1^0 \in \Pi_1S, x_2^0\in \Pi_2S, \ldots, x_{n-1}^0$ $\in
\Pi_{n-1}S$ such that for all $f\hbox{\rm :}\ S \rightarrow {\mathbb
{C}}$ the equation
\begin{align}
u_1(x_1) + u_2(x_2) + \cdots + u_n(x_n) &= f(x_1, x_2,\ldots,
x_n),\nonumber\\
&\quad\ (x_1, x_2, \ldots, x_n) \in S,
\end{align}

$\left.\right.$\vspace{-1.5pc}

\pagebreak

\noindent subject to
\begin{equation}
u_1(x_1^0) = 0,\quad u_2(x_2^0)= 0, \ldots, u_{n-1}(x_{n-1}^0) = 0
\end{equation}
admits a unique solution.  Then $S$ is full.
\end{therr}

\begin{proof}
Before we proceed with the proof we remark that the solution is unique
only in the sense that the functions $u_i|_{\Pi_iS}, 1\leq i\leq n$,
are uniquely determined and how any of the $u_i$ defined outside $X_i
- \Pi_iS$ is immaterial.
\end{proof}

Clearly $S$ is a good set since for all $f\hbox{:}\ S \rightarrow {\mathbb{C}}$,
(1) admits a solution by assumption. We show that under the given
hypothesis $S$ is full. If $S$ is not full, then there exists $\vec a =
(a_1, a_2, \ldots, a_n)$ in the Cartesian product of $\Pi_i S, 1\leq i
\leq n$, such that $S' = S \cup \{\vec a\}$ is a good set. Consider the
function $f$ on $S'$ which vanishes everywhere on $S$ and equals one at
$\vec a$. Let $U_i$, $1\leq i \leq n$, be a solution of
\begin{align}
u_1(x_1) + u_2(x_2) + \cdots + u_n(x_n) &= f(x_1,x_2,\ldots,
x_n),\nonumber\\
&\quad\ (x_1,x_2, \ldots, x_n) \in S'.
\end{align}
Then the system of functions
\begin{equation*}
V_i = U_i - U_i(x_i^0),\quad 1 \leq i \leq n-1,\quad V_n = U_n + \sum_{i
= 1}^{n-1}U_i(x_i^0),
\end{equation*}
is also a solution of (3). In particular, this system, when restricted
to $S$, is the unique solution of (1) subject to (2) for the identically
null function on $S$ (observe that $f$ vanishes on $S$), whence we have
$V_i(x_i) = 0, x_i \in \Pi_i S, 1\leq i \leq n$. Since $a_i\in \Pi_iS,
1\leq i\leq n$ we see that $\sum_{i=1}^nV_i(a_i) = 0 \neq 1$, which is a
contradiction. So $S$ is full, and the theorem is\break proved.

\begin{therr}[\!] Let $S\subset \Omega$ be full and fix $x_i^0\in \Pi_i
S, 1\leq i\leq n-1$. Then the equation
\begin{equation}
u_1(x_1) + u_2(x_2) + \cdots + u_n(x_n) = 0,\quad (x_1,x_2, \ldots, x_n)
\in S,
\end{equation}
subject to
\begin{equation}
u_1(x_1^0) = 0,\quad u_2(x_2^0)= 0, \ldots, u_{n-1}(x_{n -1}^0) = 0
\end{equation}
admits a unique solution which is necessarily the trivial solution
$U_i(x_i) = 0, x_i \in \Pi_i S$, $1\leq i \leq n$.
\end{therr}

\begin{proof} We have to show that any solution $U_1, U_2, \ldots, U_n$
of (4) subject to (5) is necessarily the trivial solution $ U_i(x_i) =
0, x_i \in \Pi_iS, 1\leq i\leq n$. If not there is a non-trivial
solution $V_i, 1\leq i\leq n$, of (4) along with (5), which means that
there exists an element ${\vec{a}} = (a_1,a_2,\ldots, a_n)\in S$ with at
least one (hence two or more) $V_1(a_1), V_2(a_2), \ldots, V_n(a_n)$
non-zero and $\sum_{i = 1}^nV_i(a_i) = 0$.

Without loss of generality assume that $V_n(a_n) \neq 0$. Since
$\sum_{i=1}^{n-1} V_i(x^0_i) + V_n(a_n) \neq 0$, $ \vec b
=(x_1^0,x_2^0,\ldots,x_{n-1}^0, a_n) \notin S$. Also ${\vec {b}}$ is in
the Cartesian product of $\Pi_iS, 1\leq i \leq n$. Consider $S' =
S\cup \{\vec b\}$. Note that $S'$ and $S$ have the same canonical
projections on the coordinate spaces. We show that $S'$ is a good set,
conflicting with the fact that $S$ is full. To this end let $f\hbox{:}\ S'
\rightarrow {\mathbb {C}}$ be given. Write $ f({\vec {b}}) = k $ and let
$W_1, W_2,\ldots, W_n$ be a solu-\break tion of
\begin{align*}
u_1(x_1) + u_2(x_2) + \cdots + u_n(x_n) &= f(x_1, x_2, \ldots, x_n),\\
&\quad\ (x_1, x_2, \ldots, x_n) \in S,
\end{align*}
subject to $u_1(x_1^0) = 0, u_2(x_2^0)= 0, \ldots, u_{n-1}(x_{n-1}^0) =
0$ which exists since $S$ is good. Write $c = \frac{k - W_n (a_n)}{V_n
(a_n)}$. Then 
\begin{equation*}
R_1 = W_1 + cV_1,\quad R_2 = W_2 + cV_2, \ldots, 
R_n = W_n + cV_n
\end{equation*}
is a solution of 
\begin{align*}
u_1(x_1) + u_2(x_2) + \cdots + u_n(x_n) &= f(x_1, x_2, \ldots, x_n),\\
&\quad\ (x_1, x_2, \ldots, x_n) \in S',
\end{align*}
which shows that $S'$ is a good set, a contradiction. The theorem is
proved.
\end{proof}

We can combine Theorems~1 and 2 as:

\begin{therr}[\!] A good set $S\subset \Omega$ is full if and only if
for any choice of $ x_i^0 \in \Pi_iS, 1\leq i \leq n-1$, the
equation
\begin{equation*}
u_1(x_1) + u_2(x_2) + \cdots + u_n(x_n) = 0,\quad (x_1, x_2, \ldots, x_n)
\in S,
\end{equation*}
subject to $u_1(x_1^0) = 0, u_2(x_2^0)= 0, \ldots, u_{n-1}(x_{n-1}^0) =
0$ has a unique solution{\rm ,} namely the trivial solution.
\end{therr}

Note that in Theorem~3 the words `any choice' can be replaced by `some
choice'.

\setcounter{theorr}{0}
\begin{coro}$\left.\right.$\vspace{.5pc} 

\noindent Let $S\subset \Omega$ be given. Then $S$ is full if and only
if for any choice of $ x_i^0 \in \Pi_iS, ~1\leq i \leq n-1${\rm ,} for
all complex valued functions $f$ on $S${\rm ,} for all complex $c_1,
c_2, \ldots, c_{n-1}${\rm ,} the equation
\begin{align*}
u_1(x_1) + u_2(x_2) + \cdots + u_n(x_n) &= f(x_1, x_2, \ldots, x_n),\\
&\quad\ (x_1,x_2, \ldots, x_n) \in S,
\end{align*}
subject to $u_1(x_1^0) = c_1, u_2(x_2^0) = c_2, \ldots,
u_{n-1}(x_{n-1}^0) = c_{n-1}$ has a unique solution.
\end{coro}

\setcounter{theorr}{0}
\begin{rem}{\rm There is nothing special about the choice of the first $n-1$
coordinates $x_1^0, x_2^0, \ldots, x_{n-1}^0$ in the sense that we could
just as well have chosen any $n - 1$ coordinates $x_i \in \Pi_iS, i \neq
i_0$, and modified the `boundary condition' accordingly.}
\end{rem}

\setcounter{theorr}{1}
\begin{coro}$\left.\right.$\vspace{.5pc}

\noindent Let $S\subset \Omega$ be full and let $U_1, U_2, \ldots, U_n$
be a solution of
\begin{equation*}
u_1(x_1) + u_2(x_2)+\cdots +u_n(x_n) = 0,\quad (x_1, x_2, \ldots, x_n)
\in S
\end{equation*}
then $U_1, U_2, \ldots, U_n$ are constant on $\Pi_1S, \Pi_2S, \ldots,
\Pi_nS$ respectively with the sum of the constants equal to zero.
\end{coro}

A corollary of the above corollary is:

\begin{coro}$\left.\right.$\vspace{.5pc} 

\noindent Let $S\subset \Omega$ be full. Let $\{1,2,\ldots,n\} = A\cup
B${\rm ,} $A\cap B = \emptyset$. Let $U_1, U_2, \ldots, U_n$ be a
solution of
\begin{equation*}
u_1(x_1) + u_2(x_2) + \cdots + u_n(x_n) = 0,\quad (x_1, x_2, \ldots,
x_n) \in S,
\end{equation*}
subject to $u_i(x_i^0) =0, i \in A$. Then $U_i(x_i) = 0 $ for all $x_i
\in \Pi_iS, i\in A${\rm ,} while if $c_j = U_j(x_j), x_j \in \Pi_jS${\rm
,}\ \,for $j\in B${\rm ,} then $\sum_{j\in B}c_j = 0$. More generally{\rm
,} if $U_1, U_2, \ldots, U_n$ and $V_1, V_2, \ldots, V_n$ are two
solutions of
\begin{align*}
u_1(x_1) + u_2(x_2) + \cdots + u_n(x_n) &= f(x_1, x_2, \ldots, x_n),\\
&\quad\ (x_1, x_2, \ldots, x_n) \in S,
\end{align*}
subject to $u_i(x_i^0) = c_i, i \in A${\rm ,} then $U_i(x_i) = V_i(x_i)
$ for all $x_i \in \Pi_iS, i\in A${\rm ,} while $ U_j(x_j)$ $-V_j(x_j)$ is
constant on $\Pi_jS$ for $j \in B${\rm ,} and if this constant be
$d_j${\rm ,} then{\rm ,} $\sum_{j\in B}d_j = 0$.
\end{coro}

If $A$ and $B$ are two subsets of $\Omega$ and if $\Pi_iA\cap \Pi_iB
\neq \emptyset$ then we say that $A$ and $B$ have a common coordinate
of the $i$th kind.

\begin{definit}$\left.\right.$\vspace{.5pc}

\noindent {\rm Two subsets $S_1,S_2$ of $\Omega$ are said to have a
common coordinate if at least one of the $n$ intersections $\Pi_iS_1\cap
\Pi_iS_2$, $1\leq i \leq n$, is non-empty. We say that $S_1,S_2$ have
$k$ distinct coordinates in common or $k$ different kinds of coordinates
in common, if at least $k$ of the above $n$ intersections are
non-empty.}
\end{definit}\vspace{.5pc}

We now make a series of set theoretic observations on full sets:
\begin{enumerate}
\renewcommand\labelenumi{(\arabic{enumi})}
\leftskip .1pc
\item If $S_1$ and $S_2$ are full, $S_1\cup S_2$ is good, and $S_1$ and
$S_2$ have $n-1$ distinct coordinates in common, then $S_1\cup S_2$ is
full.

\item If $S_\alpha$, $\alpha \in I$, is an indexed family of full sets
such that (i) $\cup_{\alpha \in I}S_\alpha$ is a good set, (ii) given
$S_\alpha, S_\beta$ in the family, there exist $S_1,S_2,\ldots, S_n$ in
the family such that $S_1 = S_\alpha, S_n = S_\beta$, and for each $i$,
$1\leq i\leq n -1$, $S_i$ and $S_{i+1}$ have $n-1$ distinct coordinates
in common, then $\cup_{\alpha\in I}S_\alpha$ is a full set.

\item The union of a totally ordered (under inclusion) family of full
sets is a full set.

\item If $S$ is a good set and $\vec x \in S$, then the union of all
full subsets of $S$ containing $\vec x$ is a full set. It is the largest
full subset of $S$ containing $\vec x$. We denote it by $F(\vec x)$ or
$F(x_1,x_2, \ldots, x_n)$.

\item If $\vec y \in F(\vec x)$ then $F(\vec y) = F(\vec x)$, for then
$F({\vec{x}})$ and $F({\vec{y}})$ have $n$ coordinates in common all of
different kind.

\item For $\vec x,\vec y \in S$, either $F(\vec x) = F(\vec y)$ or
$F(\vec x)\cap F(\vec y) = \emptyset$. Further, since $\vec x$ is always
an element of $F(\vec x)$, we see that the collection $F(\vec x)$, $\vec
x\in S$, is a partition of $S$, which we call the partition of $S$ into
full components and call $F(\vec x)$ a full component\break of $S$.

\item Two distinct full components of a good set $S$ can have at most
$n-2$ different kinds of coordinates in common.\vspace{-1pc}
\end{enumerate}

\section{Boundary set and its existence}

As a matter of convenience we will assume henceforth that the sets $X_i,
1\leq i\leq n$, are pairwise disjoint.

\begin{definit}$\left.\right.$\vspace{.5pc}

\noindent {\rm Let $S\subset\Omega$ be a good set. A subset $B \subset
\cup_{i=1}^n\Pi_iS$ is said to be a boundary set for $S$ if for any
complex valued function $U$ on $B$ and for any $f\hbox{:}\ S\rightarrow
{\mathbb{C}}$ the equation
\begin{align*}
u_1(x_1) + u_2(x_2) + \cdots + u_n(x_n) &= f(x_1, x_2, \ldots, x_n),\\
&\quad\ (x_1, x_2, \ldots, x_n) \in S,
\end{align*}
subject to $u_i|_{B\cap{\Pi_iS}} = U|_{B\cap {\Pi_iS}}, 1\leq i \leq n$,
admits a unique solution.}
\end{definit}

\subsection*{\it Examples}\vspace{-.7pc}

\begin{enumerate}
\renewcommand\labelenumi{(\arabic{enumi})}
\leftskip .1pc
\item If $S$ is full then any set of $n-1$ different kinds of
coordinates of $S$ is a boundary set\break of $S$.

\item If no two distinct full components of $S$ have a common coordinate
then $B = \cup_{i=1}^{n-1}\Pi_iC$ is a boundary set for $S$, where $C$
is any set which intersects each full component in exactly one point.

\item In case $n=2$, the full components of $S$ are the same as the
equivalence classes of the relation $L$ defined in Example~3 of \S1, the
so-called linked components in the terminology of \cite{3}. In this case
two distinct linked components have disjoint canonical projections and
the boundary set is easily described as $\Pi_1C$ where $C$ is a
cross-section of the linked components. The difficulty for the higher
dimensional case $(n > 2)$ results from the fact that two distinct full
components can admit common coordinates (although no more than $n-2$ of
distinct kind).
\end{enumerate}

\setcounter{theorr}{0}
\begin{propo}$\left.\right.$\vspace{.5pc}

\noindent Let $S\subset \Omega$ be a good set which is not full. Assume
that there exists a full set $F, S\subset F${\rm ,} such that $F -S$ is
full{\rm ,} $\Pi_iS = \Pi_iF, 1\leq i\leq n$. Then $B =
\cup_{i=1}^n\Pi_i(F-S)$ is a boundary set for $S$.
\end{propo}

\begin{proof} Let $U_i, 1\leq i \leq n$, be any complex valued functions
on $\Pi_i(F-S), 1\leq i \leq n$, respectively. Let $f\hbox{:}\ S \rightarrow
{\mathbb{C}} $ be arbitrary and extend $f$ to all of $F$ by setting
\begin{align*}
f(x_1, x_2, \ldots, x_n) &= U_1(x_1) + U_2(x_2) + \cdots + U_n(x_n),\\
&\quad\ (x_1, x_2, \ldots, x_n) \in F-S.
\end{align*}
Fix $(x_1^0, x_2^0, \ldots, x_n^0) \in F-S$. Since $F$ is full, the
equation
\begin{align}
u_1(x_1) + u_2(x_2) + \cdots + u_n(x_n) &= f(x_1, x_2, \ldots,
x_n),\nonumber\\
&\quad\ (x_1, x_2, \ldots, x_n) \in F,
\end{align}
subject to
\begin{equation}
u_1(x_1^0) = U_1(x_1^0),\quad u_2(x_2^0) = U_2(x_2^0), \ldots,
u_{n-1}(x_{n-1}^0) = U_{n-1}(x_{n-1}^0),
\end{equation}
admits a unique solution, say, $V_1, V_2, \ldots, V_n$. Since $U_i,
1\leq i\leq n$, is already a solution of
\begin{align*}
u_1(x_1) + u_2(x_2) + \cdots + u_n(x_n) &= f(x_1, x_2,\ldots, x_n),\\
&\quad\ (x_1, x_2, \ldots, x_n) \in (F-S),
\end{align*}
subject to $u_1(x_1^0) = U_1(x_1^0), u_2(x_2^0) = U_2(x_2^0), \ldots,
u_{n-1}(x_{n-1}^0) = U_{n-1}(x_{n-1}^0)$, and since $F-S$ is full, this
solution is unique and we see that
\begin{equation*}
V_i|_{\Pi_i(F-S)} = U_i,\quad 1\leq i \leq n.
\end{equation*}
We now show that $V_i|_{\Pi_iS}, 1\leq i \leq n$, is the unique
solution of
\begin{equation}
u_1(x_1) + u_2(x_2) + \cdots + u_n(x_n) = f(x_1, x_2, \ldots, x_n),
~(x_1, x_2, \ldots, x_n) \in S,
\end{equation}
subject to
\begin{equation}
u_i|_{\Pi_i(F - S)} = U_i,\quad 1\leq i \leq n.
\end{equation}
For, if $W_i, 1\leq i\leq n$, is another solution of (8) subject to (9)
distinct from $V_i, 1\leq i \leq n$, then $W_i, 1\leq i\leq n$, is also
a solution of (6) subject to (7), which is a contradiction, since this
system has a unique solution as $F$ is full. The theorem follows.
\end{proof}

We see from this theorem that to prove the existence of a boundary set
$B$ for a non-full good set $S\subset \Omega$, it is enough to prove the
existence of a full set $F$ containing $S$, having the same canonical
projections as $S$, and such that $F-S$ is also full. We have:

\setcounter{theorr}{3}
\begin{therr}[\!]
Let $S \subset \Omega$ be a good set which is not full. Then there
exists a full set $F$ containing $S$ such that {\rm (i)} $\Pi_i(S) =
\Pi_iF, 1\leq i \leq n${\rm ,} {\rm (ii)} $F -S$ is full.
\end{therr}

\begin{proof}
Since $S$ is not full there exists a ${\vec{b}} = (b_1, b_2, \ldots,
b_n)\notin S$, $b_i \in \Pi_iS, 1\leq i \leq n$, such that $S' = S\cup
\{{\vec{b}}\}$ is good. Note that $S' - S$ is a singleton, so a full
set, and the canonical projections of $S$ and $S'$ on coordinate spaces
agree.

Let ${\cal{F}}$ be the collection of good supersets $F$ of $S$ such that
\begin{enumerate}
\renewcommand\labelenumi{(\roman{enumi})}
\leftskip .1pc
\item $\Pi_i(F) = \Pi_iS, 1\leq i \leq n$,

\item $F - S$ is full.
\end{enumerate}\vspace{-.7pc}

Note that ${\cal{F}}$ is non-empty since $S'$ belongs to it. We
partially order ${\cal{F}}$ under inclusion and observe that every chain
in ${\cal{F}}$ admits an upper bound, namely the union of the members of
the chain. By Zorn's lemma ${\cal{F}}$ admits a maximal set. Let $F$ be
one such maximal set. Clearly $F$ satisfies conclusions (i) and (ii) of
the theorem since $F$ is in ${\cal {F}}$. What remains to be proved is
that $F$ is full. If $F$ is not full, there exist a non-trivial solution
$U_1, U_2, \ldots, U_n$ of
\begin{equation*}
u_1(x_1) + u_2(x_2) + \cdots + u_n(x_n) = 0,\quad (x_1, x_2, \ldots,
x_n) \in F,
\end{equation*}
subject to $u_1(x_1^0) = 0, u_2(x_2^0) = 0, \ldots, u_{n-1}(x_{n-1}^0) =
0\ ({\rm{hence ~also}} ~U_n (x_n^0) = 0)$, where $(x_1^0,x_2^0,\ldots,
x_n^0) \in (F - S)$ is fixed. Let $\vec a = (a_1, a_2, \ldots, a_n)$ be
a point in $F$ such that for some $i$, $U_i(a_i)\neq 0$. Such a point
exists since $U_i$'s form a non-trivial solution. Moreover, $\vec a$
cannot be in $F - S$ since $F-S$ is full and there the solution is the
trivial solution. Assume without loss of generality that $U_1(a_1) \neq
0$. Consider the point $ \vec b = (a_1, x_2^0, \ldots, x_n^0)$, which is
not in $F$. The set $H = F\cup \{\vec b\}$ can be shown to be a good set
as in Theorem~2. Also $\Pi_iH = \Pi_iF = \Pi_iS$ for $1 \leq i\leq n$.
Finally $H-S$ is a full set for if $V_1, V_2, \ldots, V_n$ is a solution
of
\begin{equation*}
u_1(x_1) + u_2(x_2) + \cdots + u_n(x_n) = 0,\quad (x_1, x_2, \ldots,
x_n) \in H -S,
\end{equation*}
subject to $u_1(x_1^0) = 0, u_2(x_2^0) = 0, \ldots, u_{n-1}(x_{n-1}^0) =
0\ ({\rm {hence ~also}}\ \ U_n(x_n^0) = 0)$, then it is also a solution
on $F-S$, and since $F-S$ is full, the $V_i$'s are identically zero on
$\Pi_i(F -S), 1\leq i\leq n$. Clearly, since $V_i(x_i^0)=0$, for $2\leq
i\leq n$, we see that $V_1(a_1) = 0$, so that $V_i, 1\leq i \leq n$ is a
trivial solution on $H -S$ as well, so that $H - S$ is a full set. Thus
$H$ belongs to the family ${\cal{F}}$, and is strictly bigger than the
maximal $F$, a contradiction. So $F$ is a full set. The theorem is
proved.
\end{proof}

\subsection*{\it Remarks}\vspace{-.7pc}

\begin{enumerate}
\renewcommand\labelenumi{(\arabic{enumi})}
\leftskip .1pc
\item Let $B$ be a boundary of a good set $S$ which is not full and
assume that for each $i$, $B_i = \Pi_i B\cap X_i \neq \emptyset$. Such a
boundary always exists for a non-full good set $S$. For each $i$ choose
a $b_i \in B_i$, and let $R = \cup_{i=1}^n\{b_1\}\times\{b_2\}\times
\cdots \times\{b_{i-1}\}\times\{B_i\}\times \{b_{i+1}\}\times \cdots
\times\{b_n\}$. It is easy to verify that (1) $R$ is a full set, (2) $F
= S\cup R$ is a full set with $\Pi_iF = \Pi_iS, 1\leq i \leq n$. We will
denote the full set $F$ thus obtained by $F(S,B)$ and call $F(S,B)$ a
full set associated to $S$.

\item If $B$ is a boundary of $S$ then no proper subset of $B$ can be a
boundary of $S$, hence also no proper superset of $B$ can be a boundary
of $S$.

\item Corollary 3 suggests an equivalence relation $E_i$ on $\Pi_iS$,
which is related to the notion of boundary.
\end{enumerate}

Write $xE_iy$ if there exists a finite sequence $R_1,R_2,\ldots,R_k$ of
related components such that $x\in R_1,\ y\in R_k$ and
$\Pi_iR_j\cap\Pi_iR_{j+1}\neq \emptyset$ for $1 \leq j \leq k-1$. We
call the equivalence classes of $E_i$ the $E_i$-{\it components} of
$\Pi_iS$. It is clear that a boundary $B$ of $S$ can intersect an
$E_i$-component of $\Pi_iS$ in at most one point.

We will write $E$ for the equivalence relation on $\cup_{i=1}^n\Pi_iS$
which, for each $i$, agrees with $E_i$ on $\Pi_iS$. For any set
$A\subset \Pi_iS$ we write $s_i(A)$ for the saturation of $A$ with
respect to the equivalence relation $E_i$, the symbol $s(A)$ denotes the
saturation of $A$ with respect to the equivalence relation $E$.

In a discussion with Gowri Navada it emerged that the boundary of a good
set $S$ can be described in terms of the equivalence relations $E_i, i=
1, 2, \ldots,n$ as follows:

Let $S$ be a good set and $R_\alpha, \alpha \in I$ be its related
components. Let $J_1, J_2, \ldots, J_n$ denote the set of equivalence
classes of $E_1, E_2,\ldots, E_n$. Let $C$ be a set which meets each
$R_\alpha $ in exactly one point and let $(x^\alpha_1, x^\alpha_2,
\ldots, x^\alpha_n)$ denote this point in $R_\alpha\cap C$. Note that
$J_i = \{s_i ( x^\alpha_i ):\ \alpha \in I\}$.

Let $U_1,U_2, \ldots, U_n$ be a solution for the zero function on $S$.
Then $U_i$ is a constant on $s_ i (x^\alpha_ i) $ and if we denote this
constant by $a_i^\alpha$, then we can identify $a_i^\alpha$ with $s_i
(x_i^\alpha)$ and think of $s_i (x_i^\alpha)$ as a variable, which
satisfies the relations $\sum_{i=1}^na_i^\alpha = 0$. The set of formal
finite linear combinations (with complex coefficients) of $s_i
(x_i^\alpha)$'s, which is the same as the finite linear combinations of
$a_i^\alpha$'s is a linear space for which $a^\alpha_i$, $i = 1,
2,\ldots, n$, form a generator but not a basis in view of the relations
$\sum_{i = 1}^na_i^\alpha = 0$. But we can choose a basis from among the
generators and if $B$ denotes such a basis, a selection of one point
from each element of $B$ forms a boundary of $S$. This way of getting
the boundary is more in line with the case $n = 2$, since $C$ plays a
role here.

Let $D$ be a set which meets each element of $B$ in exactly one point.
We show that $D$ forms a boundary for $S$. Let $U$ be any function on
$D$ and $U_i$ the restriction of $U$ to $D\cap \Pi_iS$. We show that
zero function on $S$ has a unique solution $U_1, U_2, \ldots, U_n$ which
agrees with $U_i$ on $D\cap \Pi_iS$. If $x_i \in D\cap \Pi_iS$ and $y_i
\in s_i(x_i)$ then define $U_i(y_i) = U_i(x_i)$. We may view $U$ as
defined on $B$. Let $z = z_j \in \Pi_jS$ and suppose $s_j(z_j) = \sum
c_k b_k$ where $b_k\in B$. We define $U_j(z_j) = \sum c_kU_k(b_k)$. This
extends $U$ to all of $\cup \Pi_iS$.

Now formal relation, $\sum_{i=1}^n s_i(x_i^\alpha) = 0$, implies that
when we replace $s_i(x_i^\alpha)$ by a finite linear combination of
$b_j$'s, that sum of the coefficients vanishes, this in turn implies
that $\sum U_i(x_i^\alpha) = 0$, and this solution of the zero function
is unique subject to the prescribed boundary values.

\section{Relation, paths and geodesics}

\begin{definit}$\left.\right.$\vspace{.5pc}

\noindent {\rm Two points $\vec x, \vec y$ in a good set $S$ are said to
be {\it related} if there exists a finite subset of $S$ which is full
and contains both $\vec x$ and $\vec y$. If $\vec x$ and $\vec y$ are
related then we write $\vec x R\vec y$.}
\end{definit}\vspace{.5pc}

The relation $R$ is obviously symmetric and reflexive. It is transitive
in view of observation 1 about full sets, so that $R$ is an equivalence
relation, whose equivalence classes we call the $R$-components of $S$.
Note that $R$-components of $S$ are full subsets of $S$. However we do
not know if $R$-components are the same as full components. Gowri Navada
\cite{15} has shown that if $S$ has finitely many related components
then these components are also the full components.

\begin{definit}$\left.\right.$\vspace{.5pc}

\noindent {\rm Let $\vec x, \vec y$ be two related points of a good set
$S$. Any finite full set $F\subset S$ containing both $\vec x$ and $\vec
y$ is called a path joining $\vec x$ and $\vec y$. Any path joining
$\vec x$ and $\vec y$ of the smallest cardinality is called a geodesic.
Cardinality of a path joining $\vec x$ and $\vec y$ is called the length
of the path.}
\end{definit}

\begin{lemm} {\it $A, B, A\cup B$ are full sets and $A\cap B \neq
\emptyset${\rm ,} then $A\cap B$ is full.}
\end{lemm}

\begin{proof} 
If $A\cap B = A$ or $A\cap B = B$ then there is nothing to prove since
$A$ and $B$ are full. Assume therefore that $A -B \neq \emptyset$ and
$B-A \neq \emptyset.$ Let ${\vec{ x}}^0 = (x_1^0, x_2^0, \ldots, x_n^0)$
be an element of $A\cap B$. Let $f$ be a complex valued function on
$A\cap B$. Let $U_1, U_2, \ldots, U_n$ be a solution of
\begin{equation}
u_1(x_1) + u_2(x_2) + \cdots + u_n(x_n) = f(x_1, x_2, \ldots, x_n),\quad
{\vec{x}}\in A\cap B,
\end{equation}
subject to
\begin{equation}
u_1(x_1^0) = 0,\quad u_2(x_2^0)= 0, \ldots, u_{n-1}(x_{n-1}^0) = 0.
\end{equation}

We show that this solution is unique. Recall that the uniqueness (to be
proved) of $U_i, 1\leq i \leq n$, is only with regard to its values on
the sets $\Pi_i(A\cap B), 1\leq i \leq n$. Define
\begin{align*}
g(x_1, x_2,\ldots, x_n) &= U_1(x_1) + U_2(x_2) +\cdots + U_n(x_n),\\
&\quad\ (x_1, x_2, \ldots, x_n) \in B.
\end{align*}
Define $h({\vec{x}}) = f({\vec{x}}), {\vec{x}}\in A\cap B,h({\vec{x}}) =
0, {\vec{x}} \in A - A\cap B$. Note that $h$ depends only on $f$ and not
on the $U_i$'s. Note that $g$ and $h$ agree on $A\cap B$, so we can
define a function $\phi$ on $A\cup B$ which equals $h$ on $A$ and equals
$g$ on $B$. Let $W_1, W_2, \ldots, W_n$ be a solution of
\begin{align*}
u_1(x_1) + u_2(x_2) + \cdots + u_n(x_n) &= \phi(x_1, x_2, \ldots,
x_n),\\
&\quad\ (x_1, x_2, \ldots, x_n) \in A\cup B,
\end{align*}
subject to $u_1(x_1^0) = 0, u_2(x_2^0)= 0, \ldots, u_{n-1}(x_{n-1}^0) =
0$.

This solution is unique since $A\cup B$ is full. The functions $W_i,
1\leq i \leq n$, when restricted to $\Pi_i B, 1\leq i\leq n$, form a
solution of
\begin{align*}
u_1(x_1) + u_2(x_2) + \cdots + u_n(x_n) &= g(x_1, x_2, \ldots, x_n),\\
&\quad\ (x_1, x_2, \ldots, x_n) \in  B,
\end{align*}
subject to $u_1(x_1^0) = 0, u_2(x_2^0)= 0, \ldots, u_{n-1}(x_{n-1}^0) =
0$.

Since $B$ is full, this solution is unique, and so if agrees with the already
known solution, namely $U_i$ on $\Pi_iB$, $1 \leq i \leq n$.

Now $W_i, 1\leq i \leq n$, when restricted to $\Pi_iA, 1\leq i \leq n$,
is the solution of
\begin{align}
u_1(x_1) + u_2(x_2) + \cdots + u_n(x_n) &= h(x_1, x_2, \ldots, x_n),\nonumber\\
&\quad\ (x_1, x_2, \ldots, x_n) \in A,
\end{align}
subject to
\begin{equation}
u_1(x_1^0) = 0,\quad u_2(x_2^0)= 0, \ldots, u_{n-1}(x_{n-1}^0) = 0,
\end{equation}
and this solution is unique since $A$ is full. Moreover, since $h$
depends only on $f$ and not on $U_i$'s, we see that $W_i|_{\Pi_iA},
1\leq i \leq n$, remain the same no matter what solution $U_1, U_2,
\ldots, U_n$ of (10) subject to (11) is chosen. Let $W_i|_{\Pi_i(A)} =
V_i, 1\leq i \leq n$. We have for any $ x_i \in \Pi_i(A\cap B)$
\begin{equation*}
U_i(x_i) = W_i(x_i) = V_i(x_i),\quad 1\leq i \leq n.
\end{equation*}

We see therefore that for each $i$, the original function $U_i$ defined
on $\Pi_i(A\cap B), 1\leq i \leq n$, is unique being the restriction of
the unique solution $V _i, 1\leq i \leq n$, of (12) subject to (13).
This proves the lemma.
\end{proof}

Note that we have proved that, under the hypothesis of the lemma,
$\cup_{i=1}^n\Pi_i( A\cap B)$ is a boundary of $A - (A\cap B)$, $B - 
(A\cap B)$, and also of $(A - A\cap B)\cup (B - A\cap B)$.

\begin{therr}[\!] 
If two points $\vec x$ and $\vec y$ in a good set are related{\rm ,}
then there is only one geodesic joining them.
\end{therr}

\begin{proof}
Let $k$ be the minimum of the cardinalities of the paths joining $\vec
x$ to $\vec y$, and let $A$ and $B$ be two paths of cardinality $k$
joining $\vec x$ to $\vec y$. By the lemma above we see that $A\cap B$
is a full set containing $\vec x$ and $\vec y$, hence a path joining
$\vec x$ and $\vec y$. If $A \neq B$, then $A\cap B$ will be a path of
smaller cardinality than $k$, a contradiction. This proves the
theorem.
\end{proof}

\begin{remm}
It is interesting to note that the full set $\{(1,0,1), (1,1,0), (0,1,1),
(0,0,0)\}$ has the property that any two distinct points are at a
geodesic distance four from each other, a situation which does not arise
when $n = 2$.
\end{remm}

\section{Procedure for solution}

We now discuss a procedure for obtaining a solution $U_i, 1\leq i \leq
n$, of the equation
\begin{align*}
u_1(x_1) + u_2(x_2) + \cdots + u_n(x_n) &= f(x_1, x_2, \ldots, x_n),\\
&\quad\ (x_1, x_2, \ldots, x_n) \in S,
\end{align*}
for a given function $f$ on a good set $S$.

\setcounter{theorr}{0}
\begin{case}{\rm 
Assume that any two points in $S$ are related so that $S$ is itself the
$R$-component of $S$. Let $f\hbox{:}\ S \rightarrow {\mathbb {C}}$ be
given. Fix ${\vec {x}}^0 = (x_1^0, x_2^0, \ldots, x_n^0) \in S$. Let
${\vec {y}}= (y_1, y_2, \ldots,y_n) \in S $. Set $U_1(x_1^0) = 0,
U_2(x^0_2) = 0, \ldots, U_{n-1}(x_{n-1}^0) = 0$. We would like to
obtain, $U_1(y_1)$, $U_2(y_2), \ldots, U_n(y_n)$, so that
\begin{equation*}
U_1(y_1) + U_2(y_2) + \cdots + U_n(y_n) = f(y_1, y_2, \ldots, y_n).
\end{equation*}
To this end let
\begin{equation*}
G = \{ {\vec{x}}^1, {\vec{x}}^2, \ldots, {\vec{x}}^k\},\quad
{\vec{x}}^0 = {\vec{x}}^1,\quad {\vec{y}} = {\vec{x}}^k,
\end{equation*}
be a geodesic joining ${\vec{x}}^0 $ to ${\vec{y}}$. Let
$(x_1^j, x_2^j, \ldots, x_n^j)$ denote the coordinates of
${\vec{x}}^j, 1\leq j \leq k$. Let
\begin{equation*}
A_i = \Pi_i G,\quad 1\leq i \leq n,\quad C = (\cup_{i=1}^n A_i) -
\{x_1^0, x_2^0, \ldots, x_{n-1}^0\}.
\end{equation*}
A function defined on $G \times C$ will be called $G\times C$ matrix.
Consider the $G\times C$ matrix $M$ defined by
\begin{equation*}
M({\vec{x}}^i, c ) = 1\quad {\rm{if}} ~~c \in \{x_1^i, x_2^i, \ldots,
x_n^i\}\cap C,\quad M({\vec{x}}^i, c) = 0 ~~{\rm {otherwise}}.
\end{equation*}

$\left.\right.$\vspace{-1.5pc}

\noindent To solve
\begin{equation*}
u_1(x_1^j) + u_2(x_2^j) + \cdots + u_n(x_n^j) = f(x_1^j, x_2^j, \ldots,
x_n^j), ~1\leq j \leq n,
\end{equation*}
subject to $u_1(x_1^1) = 0, u_2(x_2^1)= 0, \ldots, u_{n-1}(x_{n-1}^1) =
0$, means to solve for a function $g$ on $C$ which satisfies
$\sum_{c\in C} M ({\vec{x}}^j,c)g(c) = f({\vec{x}}^j)$.

Since the solution is known to exist and is unique (since $G$ is a full
set), we see that $C$ has the same number of points as $G$, namely $k$,
and the $k\times k$ matrix $M$ is invertible (since the solution exists
for all $f$ on $G$). Finally $U_i(y_i) = g(y_i)= g(x_i^k), 1\leq i\leq
n.$ If we write ${\bf M}$ for the system of $G\times C$ matrices where $G$
runs over the geodesics beginning at ${\vec{x}}^0$, and $C$ the
associated set as above, then we may write the solution of
\begin{align*}
u_1(x_1) + u_2(x_2) + \cdots + u_n(x_n) &= f(x_1, x_2, \ldots, x_n),\\
&\quad\ (x_1, x_2, \ldots, x_n) \in S,
\end{align*}
subject to $u_1(x_1^0) = 0, \ldots, u_{n-1}(x^0_{n-1}) = 0$, formally as
${\bf M}^{-1}f$.}
\end{case}

\begin{case}
{\rm If no two distinct related components of $S$ admit a common
coordinate, then we could repeat the above procedure in each related
component and get a solution.}
\end{case}

\begin{case}{\rm 
If there is a pair of related components of $S$ with a common
coordinate then the solution as in Case~2 will yield solutions only
on related components, but solutions on different related components may
not match on a common coordinate. We therefore make use of the boundary
and the full set associated to $S$ (see Remark~1, \S4).}
\end{case}

Let $S$ be a good set and let $B$ be the boundary of $S$, and $F =
F(S,B)$ the full set associated to $S$. If $f$ is a complex valued
function on $S$, we extend it to $F$ by setting it equal to zero on $R =
F - S$. If $F$, which is a full set, is also its own related component
then we can solve for
\begin{align*}
u_1(x_1) + u_2(x_2) + \cdots + u_n(x_n) &= f(x_1, x_2, \ldots, x_n),\\
&\quad\ (x_1, x_2, \ldots, x_n) \in F,
\end{align*}
subject to $u_1(x_1^0) = 0, u_2(x_2^0)= 0, \ldots, u_{n-1}(x_{n-1}^0) =
0$ with $(x_1^0, x_2^0, \ldots, x_n^0) \in F$, and restrict the solution
to $S$.

\section{Remarks on convergence}

Let $S$ be a good set in which any two points are related. If $f_k, k =
1, 2, \ldots$ is a sequence of functions on $S$ converging pointwise to
a function $f$ and if, for each $k$, $U_{k,i}, 1\leq i \leq n$, is a
solution of
\begin{align*}
u_1(x_1) + u_2(x_2) + \cdots + u_n(x_n) &= f_k(x_1, x_2, \ldots,
x_n),\\
&\quad\ (x_1, x_2, \ldots, x_n) \in S,
\end{align*}
then, in general the functions $U_{k,i}, k =1, 2, \ldots$ need not
converge as $k \rightarrow \infty$. However, it is clear from the above
discussion that if the solutions are required to satisfy the boundary
condition $U_{k,i}(x_i^0) = 0, 1\leq i \leq n - 1, 1\leq k <
\infty$, then for each $i$, the sequence $U_{k,i}, k = 1, 2, \ldots$
converges pointwise on the set $\Pi_iS$ to a function $U_i$ and these
$U_i, 1\leq i \leq n$ give the unique solution of
\begin{align}
u_1(x_1) + u_2(x_2) + \cdots + u_n(x_n) &= f(x_1, x_2, \ldots,
x_n),\nonumber\\
&\quad\ (x_1, x_2, \ldots, x_n) \in S,
\end{align}
subject to
\begin{equation}
u_1(x_1^0) = 0,\quad u_2(x_2^0) = 0, \ldots, u_{n-1}(x_{n-1}^0) = 0.
\end{equation}
If $f_k, k = 1, 2, \ldots$ converge uniformly to $f$ and if there is a
uniform bound, say $l$, for the lengths of geodesics in $S$, then, for
each $i$, the convergence of $U_{k,i}, k = 1, 2, \ldots$ is also uniform
assuming of course that the solutions $U_{k,i}, 1\leq i \leq n$, satisfy
for each $i$ and $k$, $U_{k,i}(x_i^0) = 0$. (Note that for a fixed $l$
there are only finitely many $l\times l$ zero-one invertible matrices,
so their norms are bounded away from zero.)

Thus, if $S$ is its own related component and geodesics are of bounded
length then for bounded $f$ the solution of (14) subject to (15)
consists of bounded $u_i, 1\leq i\leq n$. If $S$ is not a related
component but the set $F$ associated to $S$ is a related component
whose geodesics are of bounded length, then also (14) admits bounded
solution whenever $f$ is bounded. This sufficient condition for bounded
solution is more in line with the condition for two-dimensional case,
than the necessary and sufficient condition of uniform separability due
to Sternfeld \cite{18} or conditions discussed by Sproston and Strauss
\cite{16}.

\section{Descriptive set theoretic considerations}

Now let $X_1, X_2, \ldots, X_n$ be Polish spaces equipped with their
respective Borel $\sigma$-algebras. Let $\Omega = X_1\times X_2\times
\cdots\times X_n$ be equipped with the product Borel structure. Let $S
\subset \Omega$ be a good Borel set. We will show that the equivalence
relation $R$ is a Borel equivalence relation. To this end let $S^k =
S^{\{1, 2, \ldots, k\}}$ be the $k$-fold Cartesian product of $S$ with
itself. Let $({\vec{x}}^1, {\vec{x}}^2, \ldots, {\vec{x}}^k) \in S^k$,
${\vec{x}}^i = (x_1^i, x_2^i, \ldots, x_n^i), 1\leq i \leq n$,
\begin{equation*}
G = \{{\vec{x}}^1, {\vec{x}}^2, \ldots, {\vec{x}}^k\},\quad
C = \cup_{i=1}^{n-1}(\Pi_i G - \{x_i^1\}) \cup \Pi_n G.
\end{equation*}
Let $M({\vec{x}}^1, {\vec{x}}^2, \ldots, {\vec{x}}^k)$ denote the
$G\times C$ matrix (see \S6)
\begin{equation*}
M({\vec{x}}^i, c ) = 1\quad {\rm{if}} ~~c \in \{x_1^i, x_2^i, \ldots,
x_n^i\}\cap C,\quad M({\vec{x}}^i, c) = 0 ~~{\rm {otherwise}}.
\end{equation*}

$\left.\right.$\vspace{-1.5pc}

\noindent The mapping
\begin{equation*}
K\hbox{:}\ ({\vec{x}}^1, {\vec{x}}^2, \ldots, {\vec{x}}^k)\rightarrow
M(({\vec{x}}^1, {\vec{x}}^2, \ldots, {\vec{x}}^k))
\end{equation*}
is a Borel map from $S^k$ into the space of finite matrices. An element
$({\vec{x}}^1, {\vec{x}}^2, \ldots, {\vec{x}}^k)\in S^k$ is called an
ordered geodesic of length ${k}$ between ${\vec{x}}^1$ and ${\vec{x}}^k$
if $\{{\vec{x}}^1, {\vec{x}}^2, \ldots, {\vec{x}}^k\}$ is a geodesic
between ${\vec{x}}^1$ and ${\vec{x}}^k$.

For a proper subset $J$ of $\{1, 2, \ldots, k\}$, $\Pi_J$ will denote
the canonical projection of $S^k$ onto $S^J$. In the definition of $M_k$
below, $J$ runs over all proper subsets of $\{1, 2, \ldots, k\}$ which
contain $1$ and $k$.
\begin{align*}
M_k &= \{({\vec{x}}^1, {\vec{x}}^2, \ldots {\vec{x}}^k)\in S^k: \forall
J, M(\Pi_J({\vec{x}}^1, {\vec{x}}^2, \ldots, {\vec{x}}^k)) {\rm ~is ~not
~invertible}\},\\[.2pc]
L_k &= \{({\vec{x}}^1,{\vec{x}}^2,\cdots {\vec{x}}^k)\in S^k: M
({\vec{x}}^1,{\vec{x}}^2,\cdots {\vec{x}}^k) {\rm ~is ~invertible}\},\\[.2pc]
G_k &= L_k\cap M_k.
\end{align*}
We note that $G_k$ is the set of vectors in $S^k$ which are ordered
geodesics of length $k$ between its first and the last coordinates. It
is a Borel set since $M_k$ and $L_k$ are Borel sets. Since there are
$(k-2)!$ ordered geodesics between two points when the geodesic length
between them is $k$, the maps defined by (for $k = 1, 2, \ldots$)
\begin{equation*}
\phi_k ({\vec{x}}^1, {\vec{x}}^2, \ldots, {\vec{x}}^k) = ({\vec{x}}^1,
{\vec{x}}^k),\quad k \geq 2,\quad \phi_1({\vec{x}}^1) = ({\vec{x}}^1,
{\vec{x}}^1)
\end{equation*}
from $G_k$ to $S\times S$ are finite to 1 Borel maps, so that for each
$k$, $\phi_k(G_k)$ is a Borel set. The equivalence relation $R =
\cup_{k=1}^\infty\phi_k(G_k)$ is thus a Borel equivalence relation.

We mention here some observations due to S~M~Srivastava and
H~Sarbadhikari on the nature of the relations $R$ and $E_i$.

Let $S$ be compact, second countable  and good. Then\vspace{-.3pc}

\begin{enumerate}
\renewcommand\labelenumi{(\arabic{enumi})}
\leftskip .1pc
\item The decomposition $R$ of $S$ into related components as well as
the equivalence relations $E_i$ defined in terms of related components
are $\sigma$-compact.

\item If for each related component $L$ there is a positive integer
$N_L$ such that every geodesic in $L$ is of length at most $N_L$, then
$L$ is compact. Hence, in this case, there is a $G_{\delta}$
cross-section for equivalence classes of $R$.
\end{enumerate}\vspace{-.7pc}

Assume, moreover, that $N_L$ is independent of $L$. Then $R$ is compact.
Further, let $C$ be an $E_i$ equivalence class that is of
bounded type, in the sense that there is a positive integer $M_C$ such
that for every $x, y \in C$, one needs at most $M_C$ many related
components to witness that $x E_i y$. Then $C$ is compact. Hence, if
each $C$ is of bounded type, then $E_i$ equivalence classes admit a
$G_{\delta}$ cross-section. Further, if $M_C$ is independent of $C$,
then $E_i$ equivalence classes itself is compact.

It is not clear how to combine these facts with the second description
of the boundary given at the end of \S4 to give a good sufficient
condition for the existence of a Borel measurable boundary, a hypothesis
needed in the discussion that follows. Of course if there are only
countably many $R$ equivalence classes then the boundary is countable
too, hence Borel measurable.

If $S$ is a good Borel set and if $f$ a complex valued Borel function on
$S$, the question whether one can choose the functions $U_i, 1\leq i
\leq n,$ in (14) in a Borel fashion has, in general, a negative answer
\cite{6}. We discuss conditions under which an affirmative answer is
possible.

Assume now that the related components of $S$ admit a Borel
cross-section $\Gamma$. The set $R_k$ of ordered geodesics of length $k$
beginning at a point in $\Gamma$ is a Borel set since
\begin{equation*}
R_k = \{({\vec{x}}^1, {\vec{x}}^2, \ldots, {\vec{x}}^k)\in G_k:
{\vec{x}}^1 \in \Gamma\} = (\Pi_1^{-1}\Gamma) \cap G_k.
\end{equation*}
The set $C_k = \Pi_kR_k$ is the Borel set of points in $S$ which are
joined to some point in $\Gamma$ by a geodesic of length $k$. Clearly $S
= \cup_{k=1}^\infty C_k$, the union being pairwise disjoint, where
$C_1=\Gamma$.

It is clear from the procedure given for the solution of (14)
that
\begin{enumerate}
\renewcommand\labelenumi{(\arabic{enumi})}
\leftskip .1pc
\item if $f$ is a Borel function and $S$ has only one related component,
then the solution is made of Borel functions,

\item if $S$ admits a Borel measurable boundary and the full set $F$
associated to $S$ is its own related component, then the solution of
(14) is made of Borel functions whenever $f$ is Borel,

\item if no two related components of $S$ admit a common coordinate and
the related components of $S$ admit a Borel cross-section then the
solution is made of Borel functions whenever $f$ is Borel.
\end{enumerate}

\section{Simplicial measures and sums of algebras}

Let $X_1, X_2, \ldots, X_n$ be Polish spaces, and $\Omega$ their
Cartesian product equipped with the product Borel structure. A
probability measure $\mu$ on $\Omega$ is called simplicial if it is an
extreme point of the convex set of all probability measures on $\Omega$
whose one-dimensional marginals are the same as those of $\mu$. Let
$\mu_i$ denote the marginal of $\mu$ on $X_i$, $1\leq i \leq n$. A basic
theorem of Lindenstrauss \cite{9} and Douglas \cite{4} states that a
probability measure on $\Omega$ is simplicial if and only if the
collection of functions of the form
\begin{equation*}
u_1(x_1) + u_2(x_2) + \cdots + u_n(x_n),\quad u_i\in
L_1(X_i,\mu_i),\quad 1\leq i \leq n,
\end{equation*}
is dense in $L_1(\Omega, \mu)$.

A Borel set $E \subset \Omega$ is called a set of marginal uniqueness
(briefly an MU-set) if every probability measure $\mu$ supported on $E$
is an extreme point of the convex set of all probability measures on
$\Omega$ with one-dimensional marginals same as those of $\mu$. Clearly
any Borel subset of an MU-set is an MU-set and since a loop is not an
MU-set, an $MU$-set cannot contain a loop, whence an MU-set is a good
set.

If $S$ is a good Borel set in which any two points are related and there
is a uniform upper bound for the lengths of geodesics, then every
bounded Borel function on $S$ is a sum of bounded Borel functions on
$X_1, X_2, \ldots, X_n$ respectively and since bounded Borel functions
are dense in $L^1$, we see that $S$ is a set of marginal uniqueness.

More generally it can be shown, as in the case $n = 2$ (see \cite{5,6}),
that if $S$ is a good Borel set in which any two points are related and
there is a uniform upper bound for $U_1, U_2, \ldots, U_n$ which form the
solution of (14) subject to (15) for $f$ which are indicator functions
of singletons, then $S$ is an MU-set. Of course one can replace the
hypothesis on $S$ by a similar hypothesis on $F(S,B)$ and claim that $S$
is an MU-set.

Assume now that $X_1, X_2, \ldots, X_n$ are compact metric spaces. Let
$S \subset \Omega$ be a compact set with $\Pi_i S = X_i$, for $i = 1, 2,
\ldots, n$. It is easy to see, by considering annihilators, that $C(X_1)
+ C(X_2) + \cdots + C(X_n)$ is dense in $C(S)$ if and only if $S$ is a
set of marginal uniqueness. We see therefore that if any two points of
the set $F = F(S,B)$ are related, $S$ has a Borel measurable boundary
and if geodesics lengths in $F$ are bounded above then $C(X_1) + C(X_2)
+ \cdots + C(X_n)$ is dense in $C(S)$. In fact we also have
\begin{equation*}
C(X_1) + C(X_2) + \cdots + C(X_n) = C(S).
\end{equation*}
We see this as follows: Let $f \in C(S)$, and let $U_{1,k}, U_{2,k},
\ldots, U_{n,k}, k = 1, 2, \ldots$ be a sequence of continuous
functions on $X_1, X_2, \ldots, X_n$ respectively, such that $U_{1,k} +
U_{2,k} + \cdots + U_{n,k}$ converges to $f$ uniformly. Fix ${\vec{x}}^0
= (x_1^0, x_2^0, \ldots, x_n^0) \in S$. Let
\begin{equation*}
V_{i,k} = U_{i,k} - U_{i,k}(x_i^0),\quad 1\leq i \leq n-1,\quad V_{n,k}
= U_{n,k} + \sum_{j=1}^{n-1}U_{j,k}(x_j^0).
\end{equation*}

$\left.\right.$\vspace{-1.5pc}

\noindent Then $V_{i,k}, 1 \leq i \leq n$, are continuous and their sum
converges to $f$ uniformly. But since $V_{i,k}(x_i^0) = 0, 1\leq i \leq
n-1$, we see from our remarks on convergence that each sequence
$V_{i,k}, k = 1, 2, \ldots$ of continuous functions converges uniformly
to a continuous function $V_i$ on $X_i$ and that $f$ is the sum of these
functions.

\end{document}